\documentclass{article}
\title{Rigidity of Riemannian foliations with complex leaves on K\"ahler manifolds
\footnote{MSC 2000 : 53C12, 53C24, 53C28, 53C29  \newline
Keywords : K\"ahler manifold, Riemannian foliation }}
\author{Paul-Andi Nagy}
\date{\today}
\usepackage{amsfonts,amssymb}
\usepackage{minitoc_href}

\newtheorem{teo}{Theorem}[section]
\newtheorem{lema}{Lemma}[section]

\newtheorem{rema}{Remark}[section]
\newtheorem{coro}{Corollary}[section]
\newtheorem{nr}{}[section]
\begin{document}
\maketitle
\abstract{\normalsize We study Riemannian foliations with complex leaves on  K\"ahler manifolds. The tensor $T$, 
the obstruction to the foliation be totally geodesic, is interpreted as a holomorphic section of a certain 
vector bundle. This enables us to give classification results when the manifold is compact.} 
\large
\tableofcontents
\section{Introduction}
Riemannian foliations with totally geodesic leaves and in particular Riemannian submersions with totally geodesic fibers 
are now quite well understood. Many general structure results in the theorey of Riemannian submersions are known (see \cite{Besse}, chapter 9). For particular spaces - as 
spheres or complex and quaternionic projective spaces - classification results are available \cite{Escobales, Escobales1} under some 
geometric hypothesis on the fibers. In the less explored case of pseudo-Riemannian submersions similar results are known to hold under some additional 
conditions \cite{Magid, Baditoiu, Baditoiu1}. In the case of Riemannian foliations transversal geometric assumptions were used in order to obtain classification theorems 
\cite{Tondeur1}. \par 
In a complex setting, a notion of almost Hermitian 
submersions was proposed in \cite{Wat1} but it turns out that for many classes the horizontal distribution has to be integrable \cite{Wat1, 
Falcitelli}.  A less rigid situation, even in the case of a submersion, should arise from the study of Riemannian submersions from an  
almost-Hermitian manifold. The geometric condition we need here is that the fibers  (or the leaves)  to be almost complex.
This is of interest when searching geometric structures admitting a (Riemannian) twistor construction as explained in 
\cite{BerBer}. \par
In this paper we study Riemannian foliations with complex leaves on K\"ahler manifolds. The totally geodesic case was completely described 
in \cite{Nagy1}, where it is shown that under the 
simple connectivity and completeness assumptions such an object is a Riemannian product of twistor spaces 
over positive quaternionic K\"ahler manifolds, K\"ahler manifolds and homogeneous spaces belonging to three main classes (see \cite{Sal3} for basic 
quaternionic-K\"ahler geometry). Note that for the case of the complex projective space this was already known in 
\cite{Escobales}. \par
It is then natural to investigate the non-totally geodesic case. It turns out the ambient K\"ahler geometry is sufficiently strong to 
force, at least in the compact case, the foliation to be of very special type. More precisely, our main result is the following rigidity theorem.
\begin{teo}Let $(M,g,J)$ be a compact K\"ahler manifold. If $M$ carries a Riemannian foliation ${\cal{F}}$ with 
complex leaves then $M$ is locally isometric and biholomorphic with a Riemannian product $M_1 \times M_2$ of 
K\"ahler manifolds where $M_1$ carries a totally geodesic, Riemannian foliation with complex leaves and 
$M_2$ carries a Riemannian foliation with complex leaves which is transversally integrable. Moreover, 
the foliation ${\cal{F}}$ is the Riemannian product of the latter.
\end{teo}
As it is well known, the decomposition theorem of deRham ensures that at least locally one can restrict attention to holonomy irreducible Riemannian manifolds. For 
the case of the latter, theorem 1.1 gives :  
\begin{coro}
On a compact, simply connected, irreducible K\"ahler manifold any Riemannian foliation with complex leaves is either 
totally geodesic, or transversally integrable.
\end{coro}
Note that for these rigidity results no assumption on the curvature of the metric $g$ is necesssary. In a standard fashion, conditions ensuring 
total geodesicity of a given foliation are based on bounds on, say, Ricci curvature (see \cite{Tondeur2} for examples of results of this type). 
Note also that when studying holomorphic distributions on K\"ahler manifolds conditions on the metric are necessary even in the case 
of (real) codimension $2$ \cite{Lebrun1}. \par
The paper is organized as follows. In section 2 we collect some classical facts about Riemannian foliations and then specialize to the case of 
K\"ahler manifolds. We are basically starting from O'Neill's equations for the curvature tensor and use the K\"ahler structure to 
derive differential relations between the basic tensors $A$ and $T$. In section 3 we interpret the 
tensor $T$, the obstruction to the foliation to be totally geodesic as a holomorphic section of a certain vector 
bundle and use the compacity assumption in order to obtain the splitting in theorem 1.1.
\section{Preliminaries}
We start by collecting a number of basic facts about Riemannian foliations and next we will specialize to the K\"ahler case. Let 
$(M,g)$ be a Riemannian manifold and let ${\cal{F}}$ be a foliation on $M$. We denote by ${\cal{V}}$ the integrable distribution induced by 
${\cal{F}}$. Let $H$ be the orthogonal complement of ${\cal{V}}$. We assume the foliation ${\cal{F}}$ to be Riemannian, that is 
$$ {\cal{L}}_Vg(X,Y)=0$$
whenever $X, Y$ are in $H$ and $V$ belongs to ${\cal{V}}$. Let $\nabla$ be the Levi-Civita connection 
of the metric $g$. Throughout this paper we will denote by $V,W$ vector fields in ${\cal{V}}$ and by 
$X,Y,Z$ etc. vector fields in $H$.
It is easy to verify that the formula \cite{Tondeur}
$$ \overline{\nabla}_EF=(\nabla_EF_{{\cal{V}}})_{{\cal{V}}}+(\nabla_EF_H)_H$$
defines a metric connection with torsion on $M$ (here the subscript denotes orthogonal projection on the 
subspace). The main property of this connection is that it preserves the distributions ${\cal{V}}$ and $H$. 
If $T$ and $A$ are the O'Neill's tensors of the foliation then the following relations between $\nabla$ and $\overline{\nabla}$ are known 
to hold 
$$ \begin{array}{lr}
\nabla_XY=\overline{\nabla}_XY+A_XY, \ \nabla_XV=\overline{\nabla}_XV+A_XV \\
\nabla_VX=\overline{\nabla}_VX+T_VX, \ \nabla_VW=\overline{\nabla}_VW+T_VW.
\end{array}$$
For the algebraic properties of $T$and $A$ see \cite{Tondeur}. We only recall here that $A$ is skew-symmetric 
on $H$ while $T$ is symmetric on ${\cal{V}}$.  \par 
In the rest of this paper we will asssume that 
$(M,g)$ is a K\"ahler manifold of dimension $2m$, with complex structure $J$. Moreover, we suppose 
that the foliation ${\cal{F}}$ has complex leaves, that is $J {\cal{V}}={\cal{V}}$ (then 
of course, $JH=H$). As $\nabla J=0$, it follows that $\overline{\nabla}J=0$, hence we obtain 
information about the complex type of the tensors $A$ and $T$ as follows 
\begin{nr} \hfill 
$ \begin{array}{cc}
A_X(JY)=J(A_XY), & A_{JX}V=-J(A_XV)=A_X(JV) \\
T_{JV}W=J(T_VW), & T_{JV}X=-J(T_VX)=T_V(JX).
\end{array} \hfill $
\end{nr}
We also have $A_{JX}JY=-A_XY$ and $T_{JV}JW=-T_VW$. A consequence of the last identity is that the foliation ${\cal{F}}$ is harmonic, that is the mean 
curvature vector field vanishes. \par
We will use now the K\"ahler structure on $M$, together with suitable 
curvature identities to get some geometric information about the tensors $A$ and $T$. 
\begin{lema}
Let $X,Y,Z$ be in $H$ and $V,W$ in ${\cal{V}}$. Then we have : \\
(i) $(\overline{\nabla}_XA)(Y,Z)=0$ \\
(ii) $<A_XY, T_VZ>=0$ \\
(iii) $<(\overline{\nabla}_VA)(X,Y),W>=<(\overline{\nabla}_WA)(X,Y),V>$. \\
(iv) $(\overline{\nabla}_{JX}T)(V,W)=-J(\overline{\nabla}_XT)(V,W)$
\end{lema}
{\bf{Proof}} : \\
We will prove (i) and (ii) simultaneously. Let us denote by $R$ the curvature tensor of the Levi-Civita connection 
of the metric $g$. We first recall the O'Neill formula (see \cite{Tondeur})
\begin{nr} \hfill
$ \begin{array}{lr} 
R(X,Y,Z,V)=<(\overline{\nabla}_ZA)(X,Y),V>+<A_XY,T_VZ>-\\ 
<A_YZ,T_VX>-<A_ZX, T_VY>.
\end{array} \hfill $
\end{nr}
Since $(M,g)$ is K\"ahler one has $R(JX,JY, Z,V)=R(X,Y,Z,V)$. Hence by (2.1) we easily arrive at 
$<(\overline{\nabla}_XA)(Y,Z),V>+<A_XY,T_VZ>=0$. But we know that (see \cite{Tondeur}, page 52)
$$\sigma_{X,Y,Z}<(\overline{\nabla}_XA)(Y,Z),V>=\sigma_{X,Y,Z}<A_XY, T_VZ>$$
(here $\sigma$ denotes the cyclic sum) thus $\sigma_{X,Y,Z}<A_XY, T_VZ>=0$ and further 
$$R(X,Y,Z,V)=<A_XY,T_VZ>.$$ Using again the 
$J$-linearity of $R$ we get imediately (ii), hence (i) follows. \par 
To prove (iii) we use another O'Neill's formula stating that 
\begin{nr} \hfill
$ \begin{array}{lr}
R(V,W,X,Y)=  <(\overline{\nabla}_VA)(X,Y),W>-<(\overline{\nabla}_WA)(X,Y),V>+\\
\hspace{3cm}<A_XV, A_YW> -<A_XW,A_YV>-\\
\hspace{3cm}<T_VX,T_WY>+<T_VX,T_WY>.
\end{array} \hfill $
\end{nr}
The result follows now by (2.1) and the fact that $R(V,W,JX,JY)=R(V,W,X,Y)$. The identity in (iv) 
can be proven in the same way, using this time the identity 
$$\begin{array}{lr}
 R(X,V,Y,W)=<(\overline{\nabla}_XT)(V,W), Y>+<(\overline{\nabla}_VA)(X,Y), W>+\\
\hspace{3cm}<A_XV,A_YW>-<T_VX,T_WY>
\end{array}$$ 
the fact that $R(JX,JV,Y,V)=R(X,V,Y,W)$ and (iii)
$\blacksquare$ 
\begin{rema}
(i) By the the first two assertions of lemma 2.1 we obtain that \\ 
$R(X,Y,Z,V)=0$, a condition frequentely imposed when studying Riemannian foliations (see chapter 5 of 
\cite{Tondeur} and references therein). \\
(ii) By (i) and (ii) of the previous lemma it is easy to see that $H$ satisfies the Yang-Mills condition.  \\
(iii) Using (iii) of Proposition 2.1 and \cite{Tondeur}, page 52, we get the following relation between the covariant derivatives of $A$ and $T$
\begin{nr} \hfill 
$ 2<(\overline{\nabla}_VA)(X,Y), W>=<(\overline{\nabla}_YT(V,W),X>-<(\overline{\nabla}_XT(V,W), Y>.\hfill $
\end{nr}
We will make use of this equation in the next section. 
\end{rema}
Let us denote by $\overline{R}$ the curvature tensor of the connection $\overline{\nabla}$. Another result that will be needed in the 
next section is the following :  
\begin{lema}
We have : 
$$\overline{R}(X,Y)V=2[A_X,A_Y]V+Q(X,Y)V$$
for all $X,Y$ in $H$ and $V$ in ${\cal{V}}$ where we defined $Q(X,Y)V=T_{T_VY}X-T_{T_VX}Y$. \\
\end{lema}
The proof follows from the general formulas in \cite{Tondeur}, page 100, and lemma 2.1, (iii).
\section{The harmonicity of the tensor $T$}
In this section we begin the study of the tensor $T$. Our main idea is to consider $T$ as a $S^2({\cal{V}})$-valued 
$1$-form on $M$ and then use 
lemma 2.1, (iv) to study differential equations involving $T$. The analogy we have 
constanly in mind is the well known fact that on a compact K\"ahler manifold any holomorphic 
$1$-form is closed. We first develop 
some preliminary material. We refer the reader to the discussion in section 4 of \cite{Donaldson}. Although 
our geometric context is different, the guideliness principle concerning the K\"ahler identities and relations between 
various natural differential operators is the same. \par
For each $p \ge 0$ we define $S^{2}({\cal{V}}) \otimes \Omega^p(H)$ to be the space of symmetric endomorphims 
$\alpha : {\cal{V}} \times {\cal{V}} \to \Omega^p(H)$. We also define $S^2_A({\cal{V}})$ as the subspace of 
$S^2({\cal{V}})$ consisting of tensors 
which vanish on $A_XY, X,Y$ in $H$. \par
The ordinary exterior derivative $d$ does not preserve $\Lambda^{\star}(M)$ but $d_H$, the horizontal component 
of its restriction to $\Lambda^{\star}(H)$ does. The latter can be extended to 
$S^2({\cal{V}}) \otimes \Lambda^{\star}(H)$ by setting 
$$ (d_H\alpha)(V,W)(X_0, \ldots, X_p)=
\sum \limits_{i=0}^{p}(-1)^i (\overline{\nabla}_{X_i}\alpha)(V,W)(X_0, \ldots \hat{X_i}, \ldots X_p) $$
for every $\alpha$ in $S^2({\cal{V}}) \otimes \Lambda^{p}(H)$. Using lemma 2.1, (ii) it is easy to see 
that $d_H$ preserves $S^2_A({\cal{V}}) \otimes \Lambda^p(H)$. 
The fact that the almost complex structure $J$ is integrable induces 
a splitting 
$$ d_H=\partial_H+\overline{\partial}_H$$
on $S^{2}({\cal{V}}) \otimes \Lambda^{\star}(H)$ where $\partial_H : S^{2}({\cal{V}}) \otimes \Lambda^{p,q}(H) \to 
S^{2}({\cal{V}}) \otimes \Lambda^{p,q+1}(H)$ and $\overline{\partial}_H : S^{2}({\cal{V}}) \otimes \Lambda^{p,q}(H) 
\to S^{2}({\cal{V}}) \otimes \Lambda^{p+1,q}(H)$. \par
We need now a formula relating to the anticommutator of the operators $\partial_H$ and $\overline{\partial}_H$. Let $Q$ be the tensor defined at the end of section 1. 
If $s$ belongs to $S^2_A({\cal{V}})$ we define the action of $Q$ on 
$s$ to be $Q.s$, an element of $S^2_A({\cal{V}}) \times \Lambda^2(H)$ defined by 
$(Q.s)(V,W)(X,Y)=s(Q(X,Y)V,W)+s(V,Q(X,Y)W)$. Obviously, this can be extended to give 
a linear application 
$$ {\cal{P}} : S^2_A({\cal{V}}) \otimes \Lambda^p(H) \to S^2_A({\cal{V}}) \otimes \Lambda^{p+2}(H), {\cal{P}}\alpha=
Q.\alpha $$
having the property that ${\cal{P}}(s\alpha)=Q.s \wedge \alpha $ whenever $s$ is in $S^2_A({\cal{V}})$ and 
$\alpha$ belongs to $\Lambda^p(H)$. \par
\begin{lema}
The following holds on $S^2_A({\cal{V}}) \otimes \Lambda^{p,q}(H)$ : 
$$ \partial_H \overline{\partial}_H+\overline{\partial}_H  \partial_H={\cal{P}} $$
\end{lema}
{\bf{Proof}} : \\
Let us first compute $d_H^2q$ where $q$ belongs to $S^2_A({\cal{V}})$. An easy manipulation yields 
$$ (d_H^2q)(V,W)(X,Y)=(\overline{\nabla}_{X,Y}^2q)(V,W)-(\overline{\nabla}_{Y,X}^2q)(V,W).$$
Using the Ricci identity for the connection with torsion $\overline{\nabla}$ (see \cite{Besse}, page 26) we get 
$$(\overline{\nabla}_{X,Y}^2q)(V,W)-(\overline{\nabla}_{Y,X}^2q)(V,W)=q(\overline{R}(X,Y)V,W)+
q(V,\overline{R}(X,Y)W)+2(\overline{\nabla}_{A_XY}q)(V,W). 
$$
Using now lemma 2.2 and the fact that $q$ vanishes on vectors of the form $A_XY$ with $X,Y$ in $H$ we obtain 
that 
\begin{nr} \hfill 
$(d_H^2q)(V,W)(X,Y)=2(\overline{\nabla}_{A_XY}q)(V,W)+q(Q(X,Y)V,W)+q(V,Q(X,Y)W).\hfill $ 
\end{nr}
We consider now $\alpha$ in 
$S^2_A({\cal{V}}) \otimes \Lambda^{p,q}(H)$ and let $\{e^I \}$ a local basis of basic $(p,q)$-forms in 
$\Lambda^{p,q}(H)$. We write $\alpha=\sum \limits_{I}^{}q_I e^I$ with $q_I $ in $S^2_A({\cal{V}}) $. Since 
$d^2_H$ vanishes on basic forms (see \cite{Tondeur}) we get using the multiplicative properties of $d_H$ that 
$d_H^2\alpha=\sum \limits_{I}^{}d_H^2q_I \wedge e^I$. But 
$\partial_H \overline{\partial}_H \alpha+\overline{\partial}_H  \partial_H \alpha=(d_H^2 \alpha)^{p+1, q+1}=
\sum \limits_{I}^{}(d_H^2q_I)^{1,1}\wedge e^I$. But $A_{JX}JY=-A_XY$ and $Q(JX,JY)=Q(X,Y)$ hence 
$(d_H^2q_I)^{1,1}={\cal{P}}q_I$ and the proof is finished
$\blacksquare$  \\ \par
At this stage let us recall another particular feature of K\"ahler geometry, namely the K\"ahler identities. We state 
them on $\Lambda^p(H)$ as follows : 
$$ \begin{array}{cc}
[\partial_H, L^{\star}]=-i\overline{\partial}_H^{\star}, & [\overline{ \partial}_H, L^{\star}]=i \partial_H^{\star}\\

[\partial_H^{\star}, L]=-i \overline{\partial}_H, & [{\overline{\partial}}_H^{\star}, L]=i \partial_H
\end{array}$$
where $L$ is multiplication with $\omega^H$ in $\Lambda^{2}(H)$ defined by 
$\omega^H(X,JY)=<X,JY>$. Of course these are projection of the K\"ahler identities of $M$ and, furthermore, it is easy 
to see that they hold on $S^2_A({\cal{V}}) \otimes \Lambda^{\star}(H)$ too.  \par
Let us now define $\zeta$ in $S^2_A({\cal{V}}) \otimes \Lambda^{0,1}(H)$ by $\zeta=\alpha_T+iJ\alpha_T$. Then 
\begin{lema}
(i) $\overline{\nabla}_{JX}\zeta=-i\overline{\nabla}_X \zeta$ for all $X$ in $H$\\
(ii) $\overline{\partial}_H \zeta=0$ \\
(iii) $\partial_H^{\star}\zeta=\overline{\partial}_H^{\star}\zeta=0$ \\
(iv) $\partial_H^{\star} \partial_H \zeta=-i{\cal{P}}^{\star} L \zeta$.
\end{lema}
{\bf{Proof}} : \\
(i) is a straightforward  consequence of lemma 2.1, (iv), while (ii) comes immediately by (i) and the fact that 
$2\overline{\partial}_H=d_H+iJd_HJ$ on $S^2_A({\cal{V}}) \otimes \Lambda^{0,1}(H)$. \par
To prove (iv) we use (in the classical way) the K\"ahler identities and the previous lemma. We have 
$$ \begin{array}{ll}
\partial_H^{\star}\partial_H \zeta=-i\partial_H^{\star}[\overline{\partial}_H^{\star}, L]\zeta=
-i(\partial_H^{\star}\overline{\partial}_H^{\star})L\zeta=-i\{ \partial_H^{\star},\overline{\partial}_H^{\star})\}L\zeta+i 
\overline{\partial}_H^{\star})[\partial^{\star}_H, L]\zeta=\\
-i\{ \partial_H^{\star},\overline{\partial}_H^{\star})\}L\zeta+\overline{\partial}_H^{\star}\overline{\partial}_H\zeta=
-i\{ \partial_H^{\star},\overline{\partial}_H^{\star})\}L\zeta.
\end{array} $$ 
It suffices now to dualize the equation in lemma 3.1
$\blacksquare$ \\ \par
Before proceeding to the proof of the theorem 1.1 we need one more preliminary result. 
\begin{lema}
We have ${\cal{P}}\zeta=0$.
\end{lema}
{\bf{Proof}} : \\
Obviously it suffices to show that ${\cal{P}}\alpha_T=0$. But it is straighforward to see that 
$$\begin{array}{rr}
({\cal{P}}\alpha)(V,W)(X,Y,Z)=
\alpha(Q(X,Y)V,W)(Z)+\alpha(V,Q(X,Y)W)(Z)\\
- \alpha(Q(X,Z)V,W)(Y) -\alpha(V,Q(X,Z)W)(Y) \\
+\alpha(Q(Y,Z)V,W)(X)+\alpha(V,Q(Y,Z)W)(X)
\end{array}$$
whenever $\alpha$ belongs to $S^2_A({\cal{V}}) \otimes \Lambda^1(H)$.
We have : 
$$\begin{array}{lr} 
<T_VQ(X,Y)W,Z>=
-<Q(X,Y)W, T_VX>=-<T_{T_WY}X-T_{T_WX}Y, T_VZ>=\\
\hspace{3.8cm}=<X, T_{T_VZ}(T_WY)>-<Y, T_{T_WX}T_VZ>.
\end{array}$$
But $<X, T_{T_VZ}(T_WY)>=-<T_{T_VZ}X, T_WY>=<T_WT_{T_{V}Z}X, Y>$ hence 
$$<T_VQ(X,Y)W,Z>=<T_WT_{T_{V}Z}X, Y>-<T_{T_WX}T_VZ, Y>.$$
Taking the alternate sum on $X,Y,Z$ of this formula gives now easily the result
$\blacksquare$ \\ \par
Let we assume, in the rest of this section, that the manifold $M$ is compact and then prove theorem 1.1. At first, taking 
the scalar product with $\zeta$ in lemma 3.2, (iv) and integrating over $M$ we obtain by lemma 3.3 that $\partial_H \zeta=0$ and since 
$\overline{\partial}_H \zeta$ vanishes it follows that $d_H \zeta=0$ and 
further $d_H \alpha_T=0$. Using now (2.4) we obtain that $(\overline{\nabla}_VA)(X,Y)=0$ and we conclude that 
\begin{nr} \hfill
$ (\overline{\nabla}_EA)(X,Y)=0 \hfill $ 
\end{nr}
for all $E$ in $TM$. But the fact that $d_H \alpha_T=0$ still contains usefull information. We proceed as follows. \par
Using the proof of lemma 3.1 (namely formula (3.1) and expression in a local basis of basic forms), one 
obtains after a few standard manipulations : 
\begin{nr} \hfill 
$ d_H^2=2{\cal{L}}+{\cal{P}}\hfill $
\end{nr}
on $S^2_A({\cal{V}}) \otimes \Lambda^1(H)$ where 
$$\begin{array}{cc}
({\cal{L}}\alpha)(V,W)(X,Y,Z)= 
(\overline{\nabla}_{A_YZ}\alpha)(V,W)(X)-
(\overline{\nabla}_{A_XZ}\alpha)(V,W)(Y)+
(\overline{\nabla}_{A_XY}\alpha)(V,W)(Z)+\\
\alpha(V,W)(A_XA_YZ-A_YA_XZ+A_ZA_XY)
\end{array}$$
\begin{rema}
Formulas of type (3.3) can be proven for forms of any degree and the operator ${\cal{L}}$ can be given a more concise form. Since 
only the case of $1$-forms is needed for our purposes this 
presentation makes more visual subsequent computations.
\end{rema}
\begin{lema} $A_X(T_VW)=0$.
\end{lema}
{\bf{Proof}} : \\
Let us recall first the following O'Neill formula : 
$$ R(V_1,V_2,V_3,Z)=<(\overline{\nabla}_{V_2}T)(V_1,V_3),Z>-<(\overline{\nabla}_{V_1}T)(V_2,V_3),Z>.$$
Now, by lemma 2.1, (ii) and (3.2) we get $<(\overline{\nabla}_{V_1}T)(A_XY,V_3), Z>=0$ and it follows that 
$R(V_1,A_XY, V_3,Z)=<(\overline{\nabla}_{A_XY}T)(V_1,V_3), Z>$. Since $(M,g,J)$ is K\"ahler 
$R(JV_1, A_X(JY),V_3,Z)=R(V_1,A_XY,V_3,Z)$ which yields further to 
$$ (\overline{\nabla}_{A_{JX}Y}T)(V_1, V_3)=-J(\overline{\nabla}_{A_{X}Y}T)(V_1, V_3).$$
Using this and relations (2.1) for the tensor $A$ we obtain after some computations that 
$$ ({\cal{L}}\alpha_T)(V,W)(JX,JY,Z)+({\cal{L}}\alpha_T)(V,W)(X,Y,Z)=2<T_VW, A_XA_YZ-A_YA_XZ>.$$
Or the vanishing of $d_H \alpha_T=0$ and ${\cal{P}}\alpha_T$ implies that of ${\cal{P}}\alpha_T$ hence 
$$ <T_VW, A_XA_YZ-A_YA_XZ>=0$$ 
for all $X,Y,Z$ in $H$ and $V,W$ in ${\cal{V}}$. Taking in this last equation $Y=JX$ we arrive 
at $<A_X(T_VW), A_X(JZ)>=0$ and the conclusion is straightforward
$\blacksquare$ \\ \par
For each $m$ in $M$ we define ${\cal{V}}^0_m$ to be the vectorial subspace of 
${\cal{V}}_m$ spanned by $\{A_XY : X, Y \mbox{in} \ H_m\}$ and let $H_m^0$ be the linear span of 
$\{ A_XV : X \ \mbox{in} \ H_m, V \ \mbox{in} \ {\cal{V}}_m\}$. By (3.1) and using parallel transport with respect 
to the connection $\overline{\nabla}$
we see that we obtained smooth distributions 
${\cal{V}}^0$ and $H^0$ of $TM$ which are furthermore $\overline{\nabla}$-parallel. We denote by 
${\cal{V}}^1$ resp. $H^1$ the orthogonal complement of ${\cal{V}}^0$ resp. $H^0$ in ${\cal{V}}$ resp. $H$. We moreover 
define distributions $D^i={\cal{V}}^i \oplus H^i, i=0,1$ of $TM$. They are both $ \overline{\nabla}$-parallel because 
$D^0$ is and $D^0$ is orthogonal to $D^1$ (of course $TM=D^0 \oplus D^1$, an orthogonal direct sum). Moreover, 
$D^0$ is $\nabla$-parallel by lemma 2.1, (ii) and lemma 3.2, (ii). 
As by the same reasons $T$ 
resp. $A$ are vanishing on $D^1$ resp. $D^2$ the proof of the theorem 1.1 is finished by using the decomposition 
theorem of DeRham.

$\\$
$\\$
$\\$
\begin{flushright}
Paul-Andi Nagy \\
Institut de Math\'ematiques \\
rue E. Argand 11, 2007 Neuch\^atel \\
Switzerland \\
e-mail : Paul.Nagy@unine.ch
\end{flushright}
\end{document}